\documentclass[11pt]{article}
\usepackage{amsmath}
\usepackage{amsthm}
\usepackage{amssymb}
\usepackage{graphics}
\usepackage{bm}
\usepackage{setspace}
\usepackage[numbers]{natbib}
\usepackage{authblk}

\usepackage{tikz}
\usepackage{amsfonts}
\usepackage[cal=bickham,calscaled=1.05]{mathalfa}
\usepackage{eucal}
\usepackage{latexsym}
\usepackage{mathdots}
\usepackage{mathrsfs}
\usepackage{mathtools}
\usepackage{hyphenat}
\usepackage{amscd}
\usepackage{enumerate}
\usepackage[pagebackref,hypertexnames=false, colorlinks, citecolor=blue, linkcolor=red, pdfstartview=FitB]{hyperref}
\usepackage[latin1]{inputenc}
\usepackage{tikz}
\usepackage{tikz-cd}
\usetikzlibrary{positioning}
\usepackage[numbers]{natbib}

\usepackage{fullpage}
\usepackage{color}

\newcommand{\be}{\begin{equation}}
	\newcommand{\ee}{\end{equation}}
\newcommand{\ba}{\begin{eqnarray*}}
	\newcommand{\ea}{\end{eqnarray*}}
\newcommand{\bal}{\begin{align}}
	\newcommand{\eal}{\end{align}}
\newcommand{\baln}{\begin{align*}}
	\newcommand{\ealn}{\end{align*}}
\newcommand{\bi}{\begin{itemize}}
	\newcommand{\ei}{\end{itemize}}
\newcommand{\bn}{\begin{enumerate}}
	\newcommand{\en}{\end{enumerate}}
\newcommand{\bbm}{\begin{bmatrix}}
	\newcommand{\ebm}{\end{bmatrix}}
\newcommand{\bpm}{\begin{pmatrix}}
	\newcommand{\epm}{\end{pmatrix}}
\newcommand{\bsm}{\left ( \begin{smallmatrix}}
	\newcommand{\esm}{\end{smallmatrix} \right) }
\newcommand{\bp}{\begin{proof}}
	\newcommand{\ep}{\end{proof}}

\newcommand{\ga}{\ensuremath{\gamma}}

\newcommand{\la}{\ensuremath{\lambda }}

\newcommand{\Addresses}{{
		\bigskip

		Fouad~Naderi, \textsc{Department of Mathematics, University of Manitoba}\par\nopagebreak
		\textit{E-mail address:} \texttt{naderif@myumanitoba.ca}

		\vspace{1cm}
		
}}

\newtheorem{thm}{Theorem}

\newtheorem{prop}{Proposition}

\newtheorem*{thm*}{Theorem}

\theoremstyle{definition}
\newtheorem{defn}{Definition}
\newtheorem{remark}{Remark}
\newtheorem{eg}{Example}

\title{Weak* decomposition and Radon-Nikodym theorem for quantum expectations}

\author[]{Fouad Naderi}
\affil[]{\footnotesize University of Manitoba}

\date{\vspace{-1.2cm}}

\begin{document}
	\maketitle
	\begin{abstract}
		A quantum expectation is a positive linear functional of norm one on a non-commutative probability space (i.e., a C*-algebra). For a given pair of quantum expectations $\mu$ and $\lambda$ on a non-commutative probability space $A$, we propose a definition for weak* continuity and weak* singularity of $\mu$ with respect to $\lambda$. Then, using the theory of von Neumann algebras, we obtain the natural weak* continuous and weak* singular parts of $\mu$ with respect to $\lambda$. If $\lambda$ satisfies a weak tracial property known as the KMS condition, we show that our weak* decomposition coincides with the Arveson-Gheondea-Kavruk Lebesgue (AGKL) decomposition. This equivalence allows us to compute the Radon-Nikodym derivative of $\mu$ with respect to $\lambda$. We also discuss the possibility of extending our results to the positive linear functionals defined on the Cuntz-Toeplitz operator system.
	\end{abstract}
	\section{Introduction}

	By the Riesz-Markov theorem and the Gelfand representation of unital commutative C*-algebras, one might refer to any positive linear functional (PLF) on a general C*-algebra as a non-commutative measure. However, we point out  that the actual counterparts of classical measures are weights on von Neumann algebras, which are only densely defined and maybe unbounded and maybe unbounded. On the other hand,  by Voiculescu's non-commutative probability \cite{DNV-free.var}, every member of a C*-algebra is a random variable, and applying a unital PLF   against  such a random variable should reveal its expected value. So,  for a PLF on a C*-algebra, the term quantum (or non-commutative) expectation is more reasonable than a non-commutative measure. In addition,  for any PLF $\mu$ on a unital C*-algebra $A$, we have $\| \mu\|=\mu(1)$, so we can rescale $\mu$ to be a  quantum expectation on $A$. Thus,  studying PLFs on C*-algebras is equivalent to studying quantum expectations on them. Despite this subtle distinction in terminology, we can still adapt and reconstruct certain concepts from classical measure theory within the non-commutative framework. We should mention that  PLFs on C*-algebras are important for quantum operations and quantum information theory, see \cite{GK-ncld}. 
	
	Here, we consider a unital C*-algebra $A$ equipped with a positive linear functional (PLF) $\lambda$. Next, we try to compare any other PLF $\mu$ on $A$ with respect to $\la$. In particular, we propose a definition for  weak* continuity and weak* singularity of $\mu$ with respect to $\lambda$. Then, by the theory of von Neumann algebras, we obtain a natural weak* continuous and weak* singular decomposition of $\mu$ with respect to $\lambda$.  If  $\lambda$ has some kind of weak tracial property known as the KMS condition, we  show that our weak* method is equivalent to   Arveson-Gheondea-Kavruk Lebesgue (AGKL) decomposition. The  equivalence with AGKL decomposition helps us to calculate the Radon Nikodym derivative of $\mu$ with respect to $\lambda$. Clearly, one can use our approach to compare  quantum expectations on a non-commutative probability space and obtain the relative  Radon-Nikodym derivative of them. 
	\section{Preliminaries}

	We recall some fundamental results from Banach algebra theory that will be used in proofs without explicit mention.
	
	\begin{defn}
		Let $A$ be a Banach algebra, and $X$ be a Banach space. We call $X$ a {\it left Banach $A$-Module} if $X$ is an algebraic module over $A$, and if the bilinear map
		\begin{eqnarray*}
			A \times X &\longrightarrow& X \\
			(a,x) &\longrightarrow a.x
		\end{eqnarray*}
		is bounded, i.e., $||a.x|| \leq C ||a|| ||x||$.  {\it Right Banach $A$-modules} are defined similarly. $X$ is called a {\it Banach $A$-bimodule}, if it is both left and right $A$-module plus
		$$a.(x.b)=(a.x).b, \quad a,~b \in A, ~ x\in X;$$
		that is, we have certain associativity between left and right actions. The module $X$ is {\it unit linked} if $A$ is unital and $1.x=x$ for any $x \in X$. 
	\end{defn}

	\begin{remark} \label{dual_action}
		When $X$  is a left $A$-module, we can make $X^{*}$ into a right $A$-module via
		\begin{eqnarray*}
			X^* \times A &\longrightarrow& X^{*} \\
			(f,a) &\longrightarrow & f.a~, \qquad  \langle  f.a, ~ x \rangle = \langle f, ~a.x \rangle .
		\end{eqnarray*}
		In this case, we call $X^{*}$ a {\it dual} Banach $A$-module. In addition, we can make $X^{**}$ into a left $A$-module again by declaring
		\begin{eqnarray*}
			A\times X^{**} &\longrightarrow& X^{**} \\
			(a,F) &\longrightarrow & a.F~, \qquad  \langle a.F,~ f \rangle = \langle F,~ f.a \rangle .
		\end{eqnarray*}
		Similarly, we can start with a right $A$-module $X$ and proceed as above, and define left and right actions on $X^*$ and $X^{**}$ respectively.  $\blacksquare$
	\end{remark}
	One instance of the above situation is  when $A$ acts on itself, i.e., $X=A$, through its multiplication, which leads to the Arens' definition of product on the second dual of any Banach algebra.

	\begin{defn} \label{Arens-products}
		Let $A$ be Banach algebra and  $A^{**}$ be its second dual. For any $\Phi, ~ \Psi \in A^{**}$, we denote the first Arens product of $\Phi$ and $ \Psi$   by $\Phi \square \Psi $ and define it by :
		\begin{eqnarray*}
			\langle \Phi \square \Psi , ~ f \rangle = \langle \Phi, ~\Psi \cdot f \rangle, ~ f\in A^* \\
			\langle  \Psi \cdot f,~ a  \rangle=\langle \Psi ,~ f \cdot a \rangle , ~ a\in A \\
			\langle f \cdot a ,~ b \rangle=\langle  f, ~ ab \rangle , ~ b\in A.
		\end{eqnarray*}
		We also we denote the second Arens product of $\Phi$ and $ \Psi$   by $\Phi \diamond \Psi $ and define it by :
		\begin{eqnarray*}
			\langle \Phi \diamond \Psi ,~ f \rangle = \langle \Psi,~ f \cdot \Phi \rangle, ~ f\in A^* \\
			\langle  f \cdot \Phi,~ a  \rangle=\langle \Phi ,~ a \cdot f \rangle , ~ a\in A \\
			\langle a \cdot f ,~ b \rangle=\langle  f, ~ ba \rangle , ~ b\in A.
		\end{eqnarray*}
	\end{defn}
	\begin{remark} \label{Arens-regularity}
		\begin{enumerate}
			\item     Consider the natural embedding $j:A \longrightarrow A^{**}$ and put the first Arens product on $ A^{**}$. For any $a, ~ b \in A$, let $\Phi=j(a)$ and $\Psi=j(b)$. By the first Arens product,  for any $f \in A^*$ we have that 
			$$\langle j(a) \square j(b) , f \rangle =\langle  f, ab \rangle= \langle j( ab), f \rangle;$$
			hence, $j:A \longrightarrow A^{**}$ is a homomorphism if we endow $A^{**}$ with its first Arens product.
			
			\item Alternatively, we can define the first and second Arens product by the aid of the natural embedding $j:A \longrightarrow A^{**}$. Note that by the Goldstine theorem $\overline{j(A)}^{wk*}=A^{**}$, where the weak* topology is $\sigma(A^{**}, A^*)$ topology. Now, if we let $\Phi=\mathrm{w^*lim}_{\alpha} j(a_{\alpha})$ and $\Psi=\mathrm{w^*lim}_{\beta} j(b_{\beta})$, then
			\begin{eqnarray*}
				\langle \Phi \square \Psi ,~ f \rangle = \mathrm{w^*lim}_{\alpha} \mathrm{w^*lim}_{\beta} \langle j(a_{\alpha} b_{\beta}), ~ f \rangle, ~~ f\in A^* \\
				\langle \Phi \diamond \Psi , ~ f \rangle = \mathrm{w^*lim}_{\beta} \mathrm{w^*lim}_{\alpha} \langle j(a_{\alpha} b_{\beta} ), ~  f \rangle, ~~ f\in A^*
			\end{eqnarray*}
			\item Now comes the question whether the Arens first and second product on the second dual $A^{**}$ of a Banach algebra $A$ are equal.  If they coincide, we call $A$ Arens regular; otherwise, we call it Arens irregular. One can show that for a locally compact group $G$, the convolution Banach algebra $L^1(G)$ is Arens regular if and only if $G$ is finite \cite[Example 5.1.19]{Runde_amn}. This is while that any operator algebra is Arens regular \cite[Corollar 2.5.4]{B.M-OA}.
			\item Finally, we should mention that if a Banach algebra $A$ is Arens regular, then the Arens product on $A^{**}$ is separately weak* continuous \cite[Remark 2.5.3]{B.M-OA}. $\blacksquare$
		\end{enumerate}
	\end{remark}
	
	The following definition is taken from \cite{DNV-free.var} and \cite{DSW-Lectures.free.prb}.
	
	\begin{defn} \label{NC-random}
		A \textit{non-commutative probability space} is a unital operator algebra $A$ equipped with a unital bounded linear functional $\phi$. The functional $\phi$ is referred to as the \textit{non-commutative (or quantum) expectation}, and the elements of $A$ are called \textit{ non-commutative (or quantum) random variables}.  If $A$ is a C*-algebra or a von Neumann algebra, we additionally require $\phi$ to be a state or a normal state, respectively.
	\end{defn}
	
	\begin{remark} \label{norm-of-pstv-lin-fun}
		Let $A$ be a unital C*-algebra  and $B$ be a C*-sub-algebra of $A$ containing the unit element of $A$. By \cite[Corollary 3.3.4]{Murphy}, any positive linear functional on any unital C*-algebra attains its norm at the identity element.  Thus, we have the following observations: 
		\begin{enumerate}[(i)]
			\item  Suppose that $\mu: A \longrightarrow \mathbb{C}$ is  a  positive linear functional on $A$, and $\mu^{\prime}: B \longrightarrow \mathbb{C}$ is the restriction of $\mu$ to $B$. Then, $\|\mu\|=\|\mu^{\prime}\|$.
			
			\item Conversely, suppose that $\nu: B \longrightarrow \mathbb{C}$ is  a  positive linear functional on $B$, and $\tilde{\nu}: A \longrightarrow \mathbb{C}$ is a positive linear functional extending  $\nu$ to $A$. Then, $\|\tilde{\nu}\|=\|\nu\|$.
		\end{enumerate}
		
	\end{remark}

	\begin{defn}\label{C-ast-dyn-sys}
		A C*-dynamical system $(\mathscr{A}, G, \alpha)$ consists of a C*-algebra $\mathscr{A}$, a locally compact group $G$, and a $\ast$-homomorphism $\alpha: G \to \operatorname{Aut}(\mathscr{A})$ such that the action  
		$$  
		G \times \mathscr{A} \longrightarrow \mathscr{A}, \quad (t, a) \mapsto \alpha_t(a),
		$$  
		is continuous with respect to the usual topologies of $G$ and $\mathscr{A}$. 
		A covariant representation of $(\mathscr{A}, G, \alpha)$ is a pair $(\pi, \rho)$, where
		\begin{enumerate}
			\item $\pi$ is a $\ast$-representation of $\mathscr{A}$ into $\mathscr{B}\left( \mathscr{H}\right)$
			\item $\rho$ is a unitary representation of $G$ into $\mathscr{B}\left( \mathscr{H}\right)$  
			\item  The {\it ``covariant covariant ( or intertwining ) relation"}  holds, i.e.,  
			$$\pi(\alpha_t(a)) = \lambda(t) \pi(a) \lambda(t)^*.$$  
			
		\end{enumerate} 
	\end{defn}

	\begin{defn} \label{W-ast-cross-prd}
		A W*-dynamical system $(\mathscr{A}, G, \alpha)$ consists of a von Neumann algebra $\mathscr{A}$ contained in some $\mathscr{B}(\mathscr{H})$, a locally compact group $G$, and a $\ast$-homomorphism $\alpha$ from $G$ into the $\ast$-automorphisms of $\mathscr{A}$ such that the action $G \times \mathscr{A} \longrightarrow \mathscr{A}$ given by $(t, a) \mapsto \alpha_t(a)$ is continuous when $G$ is equipped with its topology and $\mathscr{A}$ is endowed with the strong operator topology (SOT). Recall that we have the natural identification  
		$$L^2(G)\otimes_2 \mathscr{H} \cong L^2(G, \mathscr{H}).$$  
		A covariant representation of the W*-dynamical system $(\mathscr{A}, G, \alpha)$ is a pair $(\pi, \lambda)$, where:  
		\begin{enumerate}
			\item $\pi$ is a $\ast$-representation of $\mathscr{A}$ into $\mathscr{B}\left(L^2(G)\otimes_2 \mathscr{H}\right)$ given by  
			$$(\pi(a)f)(s) = \alpha_s(a) f(s),$$  
			for all $a \in \mathscr{A}$, $f \in L^2(G, \mathscr{H})$, and $s \in G$. 
			\item $\lambda$ is a unitary representation of $G$ into $\mathscr{B}\left(L^2(G)\otimes_2 \mathscr{H}\right)$ defined by  
			$$(\lambda(t)f)(s) = f(t^{-1}s).$$  
			\item  The covariant ( or intertwining )  relation holds, i.e.,  
			$$\pi(\alpha_t(a)) = \lambda(t) \pi(a) \lambda(t)^*.$$  
		\end{enumerate}
		The W*-crossed product $\mathscr{A} \overline{\rtimes}_{\alpha} G$ is defined as the von Neumann algebra generated by the images of $\pi(\mathscr{A})$ and $\lambda(G)$ inside $\mathscr{B}\left(L^2(G)\otimes_2 \mathscr{H}\right)$, that is,  
		$$\mathscr{A} \overline{\rtimes}_{\alpha} G = \left[\pi(\mathscr{A}) \cup \lambda(G)\right]^{''}.$$  
	\end{defn}
	
	\section{The weak* decomposition of quantum expectations}
	The following observation is crucial for constructing our decomposition, so we state it as a remark for later use. 
	\begin{remark} \label{weak-astr-cont-absol-cont-observation}
		Suppose that $\la$ is a classical , finite, positive, Radon measure on the unit circle $\mathbb{T}$. Consider the multiplication mapping
		\begin{eqnarray} \label{obvious-rep}
			\pi : L^\infty(\mathbb{T}, \la) &\longrightarrow & \mathscr{B}(L^2(\mathbb{T}, \la))  \nonumber\\
			\pi(f) &=&M_f, \quad M_f(g)=fg,
		\end{eqnarray}
		which is an isometric unital $\ast$-isomorphism from $L^\infty(\mathbb{T}, \la)$ onto the range of $\pi$, see \cite[Example 2.5.1]{Murphy}. By \cite[Proposition 4.51]{Douglas}, the weak* topology of $L^\infty(\mathbb{T}, \la)$ coincides with the weak operator topology of $\pi(L^\infty(\mathbb{T}, \la))$. Since $\lambda$ is a Radon measure, $C(\mathbb{T})$ can be regarded as a closed subspace (even sub-C*-algebra) of $L^\infty(\mathbb{T}, \la)$, see \cite[P. 184 and Section 7.2]{Folland}. By \cite[ Corollary 4.53]{Douglas}, von-Neumann double commutant theorem \cite[Theorem 4.1.5]{Murphy}, and \cite[Corollary 4.2.8]{Murphy},  we have  
		$$\overline{\pi(C(\mathbb{T}))}^{wk*}=\pi (C(\mathbb{T}))^{\prime \prime}=\pi( L^\infty(\mathbb{T}, \la)) \cong L^\infty(\mathbb{T}, \la).$$
		Now, let  $\mu$  be an arbitrary classical, finite, positive, Radon measure on the unit circle. By  Riesz-Markov Theorem \cite[Corollary 7.18]{Folland} and \cite[Theorem III.1.1]{Takesaki_I}, $\mu$ can be regarded as a positive linear functional (PLF) on $C(\mathbb{T})$ given by integration against $\mu$, i.e.,
		\begin{eqnarray*}
			\mu : C(\mathbb{T}) &\longrightarrow & \mathbb{C}\\
			\mu(f) &=&\int f~d\mu .
		\end{eqnarray*}
		By  \cite[Theorem III.1.2]{Takesaki_I}, we see that the GNS construction of $\la$ induces the irreducible $\ast$-representation given by \ref{obvious-rep}, i.e., 
		\begin{eqnarray} \label{GNS-lambda}
			\pi_\la : C(\mathbb{T})  &\longrightarrow & \mathscr{B}(L^2(\mathbb{T}, \la))  \nonumber\\
			\pi_\la(f) &=&M_f, \quad M_f(g)=fg.
		\end{eqnarray}
		So, in fact 
		\begin{equation} \label{GNS-lambda-generated}
			\overline{\pi_\la(C(\mathbb{T}))}^{wk*}=\pi_\la (C(\mathbb{T}))^{\prime \prime}= \pi_\la( L^\infty(\mathbb{T}, \la)) \cong L^\infty(\mathbb{T}, \la).
		\end{equation}
		On the other hand, $L^\infty(\mathbb{T}, \la)$ is a von Neumann algebra with the pre-dual $L^1(\mathbb{T}, \la)$. By the general theory of von Neumann algebras \cite[Sections II.2 and III.2]{Takesaki_I}, $L^1(\mathbb{T}, \la)$ is the space of all weak* continuous linear functionals on $L^\infty(\mathbb{T}, \la)$. 
		With this preamble in mind, we would like to explore the absolute continuity of a measure $\mu$ with respect to a measure $\la$ and its relationship with weak* continuous functionals on $L^\infty(\mathbb{T}, \la)$. In particular, we see that for a positive function $f \in L^1(\mathbb{T}, \la)$, the measure or PLF $f d\la$ defined by 
		\begin{eqnarray*} 
			f d\la : C(\mathbb{T})  &\longrightarrow & \mathbb{C}\\
			f d\la (g) &=&\int gf d\la ,
		\end{eqnarray*}
		can be extended to the PLF 
		\begin{eqnarray*} 
			f d\la : L^\infty(\mathbb{T}, \la)  &\longrightarrow & \mathbb{C}\\
			f d\la (g) &=&\int gf d\la,
		\end{eqnarray*}
		which defines a weak* continuous functional on $ L^\infty(\mathbb{T}, \la)$.  Interestingly, these are all the possible weak* functionals, and by the Radon-Nikodym theorem, they are all absolutely continuous with respect to $\la$. $\blacksquare$
	\end{remark}
	The following classical observation appears to be a well-known theorem. However, as we could not locate a specific reference, and since it is essential for our weak* decomposition of PLFs, we will state it explicitly and provide a proof.

	\begin{thm}\label{wc-ac-equiv--classical-measure}
		Let $\mu$ and $\lambda$ be two positive classical Radon measures on the unit circle, and consider them as positive linear functionals on $C(\mathbb{T})$. Then,
		\begin{enumerate}[(i)]
			\item $\mu$ can be extended to a weak*-continuous functional on $ {L}^{\infty}(\lambda)$  if and only if $\mu$ is absolutely continuous with respect to $\la$.
			
			\item No extension of $\mu$  to  $ {L}^{\infty}(\lambda)$ can majorize a non-zero weak*-continuous PLF on ${L}^{\infty}(\lambda)$ if and only if $\mu$ is singular with respect to $\la$.
			\item  $\mu$ can be decomposed with respect to $\lambda$ into two parts, a weak* continuous part and a weak* singular part, i.e., $\mu=\mu_{wc} + \mu_{ws}.$
		\end{enumerate}
	\end{thm}
	\bp
	
	\begin{enumerate}[(i)]
		\item Suppose that $\mu$ and $\lambda$ are  two  classical positive  Radon measures on $\mathbb{T}$. By Riesz-Markov theorem, $\mu: C(\mathbb{T}) \longrightarrow \mathbb{C}$ is a PLF acting by integration against $\mu$. Let $\mu$  has a weak$^{*}$-continuous extension $\hat{\mu}$ to a PLF on $L^{\infty}(\mathbb{T}, \lambda)=\overline{C(\mathbb{T})}^{wk*}$. Here, wk* is the the weak* topology of $L^{\infty}(\mathbb{T}, \lambda)=(L^1(\mathbb{T}, \lambda))^*$. By the   definition of weak* topology,   any weak* continuous functional corresponds to a member of $L^1(\mathbb{T}, \lambda)$. By Example \ref{weak-astr-cont-absol-cont-observation}, we view $L^{\infty}(\mathbb{T}, \lambda)$ inside $\mathscr{B}(L^2(\lambda))$ as multiplication by bounded operators such that weak* topology of $L^{\infty}(\mathbb{T}, \lambda)$ coincides with the  weak operator topology it receives from $\mathscr{B}(L^2(\lambda))$. Since the measures in question are positive, and $\hat{\mu}$ is continuous with respect to the weak operator topology, by \cite[Theorem 4.2.6]{Murphy},  
		we can find a positive $f \in L^{1}(\mathbb{T}, \lambda)$ such that for any $g \in C(\mathbb{T}) $ 
		$$\int_{\mathbb{T}} g d\mu=\mu(g)=\langle g\sqrt{f},~ \sqrt{f} \rangle_{L^2(\lambda)} =\int_{\mathbb{T}} g f d\lambda.$$
		Hence, $\mu$ and $f d\la$ give rise to the same functional on $C(\mathbb{T})$, i.e., $d\mu = f d\lambda$. Thus, $\mu$ is absolutely continuous with respect to $\lambda$.
		
		Conversely, if $\mu$ is absolutely continuous with respect to $\lambda$, then by the Lebesgue-Radon-Nikodym theorem, there exists a positive $f \in L^{1}(\mathbb{T}, \lambda)$ such that  $d\mu = f d\lambda$. So, for any $g \in C(\mathbb{T}) $ we have 
		$$\mu(g)=\int_{\mathbb{T}} g d\mu=\int_{\mathbb{T}} g f d\lambda=\langle g\sqrt{f},~ \sqrt{f} \rangle_{L^2(\lambda)},$$
		which shows that $\mu$ is a vector state. However, the vector state $\hat{\mu}(g)=\langle g\sqrt{f},~ \sqrt{f} \rangle_{L^2(\lambda)}$ on $L^{\infty}(\mathbb{T}, \lambda)$ is  an obvious extension of $\mu$  to a weak* continuous functional on $L^{\infty}(\mathbb{T}, \lambda)$. 
		
		\item Suppose that  $\mu$ has no  extension to a PLF on $L^{\infty}(\mathbb{T}, \lambda)$ that  majorize a weak* continuous functional on $L^{\infty}(\mathbb{T}, \lambda)$. By part (i), this means that it does not majorize any $\la$-absolutely continuous measure  on $C(\mathbb{T})$. So, if we consider the Lebesgue decomposition 
		$$\mu=\mu_{ac} + \mu_{s}$$
		of $\mu$ with respect to $\lambda$, we must have $\mu_{ac}=0$, i.e., $\mu =\mu_s$. That is $\mu$ is Lebesgue singular with respect to $\lambda$. 
		
		Conversely, let $\mu$ be a Lebesgue singular measure with respect to $\lambda$. Suppose to the contrary that $\mu$ has a PLF extension $\hat{\mu}$ on $L^{\infty}(\mathbb{T}, \lambda)$ that majorizes a weak* continuous functional $\ga$. Consider the Lebesgue decomposition of $\mu$ with respect to $\la$ is $\mu=\mu_{ac} + \mu_{s}$. By part (i), $\ga$ is absolutely continuous, and by our assumption $\mu_{ac} \geq \ga \neq 0$. This violates the singularity of $\mu$ with respect to $\la$. Thus, no extension of  $\mu$ can majorize a non-zero weak* continuous functional on $L^{\infty}(\mathbb{T}, \lambda)$. 
		\item Recall that any positive linear functional $\mu$ on $C(\mathbb{T})$ has a positive extension on  $L^\infty(\la)$ with the same norm, see \cite[Theorem 3.3.8]{Murphy}. So, the decomposition in this part follows from the Lebesgue decomposition theorem, and  parts (i) , (ii).
	\end{enumerate}
	\ep
	\begin{remark}
		The decomposition $\mu=\mu_{wc}+\mu_{ws}$  on $L^\infty(\la)$ in the third part of Theorem \ref{wc-ac-equiv--classical-measure} is the same as the decomposition of a positive linear functional on a von Neumann algebra into normal and singular, see \cite[P.127]{Takesaki_I}.
		Therefore, in the classical setting of measure theory and functional analysis, the absolute continuity of a measure is equivalent to the ability of  extending the corresponding Riesz-Markov functional to a weak* continuous one. In addition, singularity of measures means the lack of this ability. So, the above theorem justifies the following definition.
	\end{remark}
	
	\begin{defn}\label{wc-classical-measure}
		Let $\mu$ and $\lambda$ be two classical, finite, positive, Radon measures on the unit circle $\mathbb{T}$ with their corresponding Riesz-Markov integration functionals on $C(\mathbb{T})$ still denoted by $\mu$ and $\la$. 
		\begin{enumerate}[(i)]
			\item  We  say that $\mu: C(\mathbb{T}) \longrightarrow \mathbb{C}$ is weak* continuous with respect to  $\lambda$ if $\mu$ has a weak* continuous extension to a PLF on $ {L}^{\infty}(\lambda)$; i.e., there is an extension $\hat{\mu}:{L}^{\infty}({\lambda})\longrightarrow \mathbb{C}$ of $\mu$ which is weak*-continuous and positive. We let $WC[\la]$ denote the set of all weak*-continuous measures with respect to $\la$.
			\item   $\mu$ is called {\it weak*-singular} with respect to  $\la $ if no extension of $\mu$ to  $ {L}^{\infty}(\lambda)$ can  majorize a non-zero weak*-continuous PLF with respect to $\lambda$.   We let $WS[\la]$ denote the set of all weak*-singular measures with respect to $\la$. 
		\end{enumerate}
	\end{defn}

	Recall that if $\tau$ is a PLF on a unital C*-algebra $A$, then its isotropic left ideal of zero length vectors is $\{a \in A: \tau(a^*a)=0\}$.
	\begin{defn}\label{Transfer-of-measures}
		Consider two PLFs $\mu$ and $\lambda$ on a unital C*-algebra $A$ with left isotropic ideals $N_\lambda \subseteq N_\mu$. Now we want to reuse the same classical notation  used in Example \ref{weak-astr-cont-absol-cont-observation}  for  the various spaces related to the GNS construction  of $\la$, e.g.,  the mapping (\ref{GNS-lambda}) and Equation (\ref{GNS-lambda-generated}).  Inspired by  \cite[Theorem III.2.1 and P. 321]{Takesaki_I}, we can generalize these notations to denote the GNS representation of $\lambda$ with the cyclic vector $\xi_\lambda = 1+N_\lambda$ by  $(\pi_\lambda, L^2(\lambda), \xi_\lambda)$. By \cite[Theorem 3.1.6]{Murphy}, $\pi_\lambda(A)$ is a C*-sub-algebra of $\mathscr{B}  (L^2(\lambda))$. Furthermore,  we denote  the von Neumann algebra generated by the C*-algebra $\pi_\lambda(A)$ inside $\mathscr{B}  (L^2(\lambda))$ by $L^\infty(\lambda)=\pi_\lambda(A)''=\overline{\pi_\lambda(A)}^{sot}$, and the predual of $L^\infty(\lambda)$ is shown by $L^1(\lambda)$. However, note that in the classical setting,   $C(\mathbb{T})$ is a closed subspace of $ L^\infty(\mathbb{T}, \la)$, and the representation  (\ref{GNS-lambda}) is isometric. In the general non-commutative setting, the PLF $\lambda$ is not typically faithful, meaning that $A$ is not isometrically embedded in $L^\infty(\lambda)$. Therefore, we need a method to transfer the PLFs $\mu$ and $\lambda$ on $A$ to PLFs on $\pi_\lambda(A)$, while ensuring that the main properties of the PLFs are preserved. This is possible by the following commuting diagrams:
		\begin{center}
			\begin{tikzcd}[row sep=large, column sep=large]
				A \arrow[r, "\pi_\lambda"] \arrow[dr, "\mu"'] 
				& \pi_\lambda(A) \arrow[d, dotted, "\mu'"] \\
				& \mathbb{C}
			\end{tikzcd}
			\quad
			\(\mu'(\pi_\lambda(a)) := \mu(a)\),
		\end{center}
		and
		\begin{center}
			\begin{tikzcd}[row sep=large, column sep=large]
				A \arrow[r, "\pi_\lambda"] \arrow[dr, "\lambda"'] 
				& \pi_\lambda(A) \arrow[d, dotted, "\lambda'"] \\
				& \mathbb{C}
			\end{tikzcd}
			\quad
			\(\lambda'(\pi_\lambda(a)) := \lambda(a)\).
		\end{center}
	\end{defn}
	
	\begin{remark}   
		Recall that in the classical setting, $\pi_\la$ is an isometric $\ast$-homomorphism, so Definition \ref{Transfer-of-measures} covers the classical case as well. However, in the general non-commutative setting, we have the short exact sequence
		\begin{center}
			\begin{tikzcd}
				0 \arrow[r] 
				& \mathrm{ker}(\pi_{\lambda}) \arrow[r, hook, "j"] 
				& A \arrow[r, "\pi_\lambda"] 
				& \pi_\lambda(A) \arrow[r] 
				& 0,
			\end{tikzcd}
		\end{center}
		and $A$ is an extension of $\pi_\lambda(A)$ by $\mathrm{ker}(\pi_{\lambda})$. Despite this, we may still assume without loss of generality that $\lambda$ is faithful, as the following argument shows. Consider the Gelfand-Naimark universal representation $\pi: A \longrightarrow \mathscr{B}(\mathscr{H})$. Recall that, by the Sherman-Takeda Theorem \cite[Lemma III.2.2]{Takesaki_I} , the enveloping von Neumann algebra $\pi(A)^{\prime\prime}$ coincides with $A^{**} = \mathscr{M}$. Since $A$ is unital and $\lambda$ is positive, we have $\|\lambda\| = \lambda(1)$ by \cite[Corollary 3.3.4]{Murphy}. By \cite[Lemma A.2.2]{B.M-OA}, $\lambda$ admits a unique linear extension $\Lambda$ to a weak* continuous functional on $\mathscr{M}$ with the same norm. Moreover, $\|\Lambda\| = \|\lambda\| = \lambda(1) = \Lambda(1)$; hence, $\Lambda$ is positive by \cite[Corollary 3.3.4]{Murphy}. By \cite[Lemma III.3.6]{Takesaki_I}, there exists a central projection $p \in \mathscr{M}$ such that $\Lambda$ is faithful on $p \mathscr{M} p$. Therefore, we may instead work with the $C^*$-algebra $pAp$ and the faithful positive functional $\lambda_p := \Lambda|_{pAp}$. In addition, by \cite[Theorem 3.3.8]{Murphy}, we can extend $\mu$ to a positive functional $M$ on $\mathscr{M}$ with the same norm, so we may work with $\mu_p := M |_{pAp}$ in place of $\mu$. This justifies that there is no loss of generality in assuming that $\lambda$ is faithful from the outset. Nevertheless, we will not make use of this assumption in the sequel.$\blacksquare$
	\end{remark}
	
	We need to show that the transferred functionals $\mu'$ and $\lambda'$ in Definition \ref{Transfer-of-measures} are well-defined. 
	
	\begin{prop}\label{wel-definedness}
		With the notations as in Definition \ref{Transfer-of-measures}, the functionals $\lambda^{\prime}$ and  $\mu^{\prime}$ are well defined and give rise to  PLFs on $\pi_\lambda (A)$. Furthermore, they have positive linear extensions to $L^\infty(\lambda)$.
	\end{prop}
	\bp
	By \cite[P. 94]{Murphy}, the GNS mapping $\pi_\la: A \longrightarrow \mathscr{B}(L^2(\la))$ defined by $\pi_\lambda (a)(b+N_\lambda)=ab+N_\lambda$ is a $\ast$-homomorphism, and by \cite[Theorem 3.1.6]{Murphy} $\pi_\la(A)$ is a C*-subalgebra of $\mathscr{B}(L^2(\la))$. Clearly,  $\pi_\la(A)$ embeds isometrically inside $L^\infty(\la)$. 
	Let $\pi_\lambda (a)=0$. Then, $\pi_\lambda (a)(b+N_\lambda)=N_\lambda$ for any $b\in A$. That is $ab \in N_\lambda$. In particular for $b=1$, we have  $a \in N_\lambda$, so $\lambda(a^*a)=0$. By the Kadison inequality \cite{Kadison-gen.CBS.inq}, we have
	$$\lambda(a)^* \lambda(a) \leq \|\lambda\| \lambda(a^* a).$$
	Thus, $\lambda(a)=0$, so $\lambda^{\prime}$ is well-defined. Since $a\in N_\lambda \subseteq N_\mu$, we have $\mu(a^*a)=0$. By another application of Kadison inequality, we see that $\mu(a)=0$, i.e., $\mu^{\prime}$ is  also well-defined. Both $\mu^{\prime}$ and  $\la^{\prime}$ are  PLFs on $\pi_{\la}(A)$ since $\pi_\lambda$ is a bounded $\ast$-homomorphism between C*-algebras. Note that $\pi_\lambda (A)$ is unital since $A$ is unital. Hence,    $\|\mu^{\prime}\|=\mu^{\prime}(1)=\mu(1)=\|\mu\|$ by \cite[Corollay 3.3.4]{Murphy}.  For the second part of the theorem we note that $\pi_\la(A)$ is a C*-sub algebra of $L^{\infty}(\lambda)$; thus, by \cite[Theorem 3.3.8]{Murphy} we can extend any PLF on $\pi_\la(A)$  to a PLF on $L^{\infty}(\lambda)$ with the same norm. 
	\ep
	\begin{remark} \label{notation-agreement}
		With Proposition \ref{wel-definedness} at our disposal, we will drop the prime notations and the homomorphism symbols in expressions like $\mu^{\prime}(\pi_\lambda(a))$ and $\lambda^{\prime}(\pi_\lambda(a)) $, and simply write $\mu(a)$ and $\lambda(a)$, respectively. Another abuse of notation is that we use the same notations $\mu$ and $\lambda$ for the extended PLFs on $L^\infty(\lambda)$ unless otherwise there is a danger of confusion. $\blacksquare$
	\end{remark}
	
	Now, Example \ref{weak-astr-cont-absol-cont-observation}, Theorem \ref{wc-ac-equiv--classical-measure}, Proposition \ref{wel-definedness}, and Remark \ref{notation-agreement} inspire us to generalize Definition \ref{wc-classical-measure} from commutative C*-algebras to the non-commutative ones. In addition, we should mention that other main sources for our inspiration were  \cite[Proposition 2.2]{DLP-ncld}, \cite[Theorem 10.1.17]{Kadison_Fundamentals_II}, \cite[Sections III.2 and III.3]{Takesaki_I},  \cite[Section III.2]{Blackadar-OA}, and some other results in Kadison and Ringrose's book, including  Proposition 10.1.20, Proposition 10.4.1, and Theorem 10.4.3, in  \cite{Kadison_Fundamentals_II}.

	\begin{defn}\label{wc-msr} 
		Consider two PLFs $\mu$ and $\lambda$ on a unital C*-algebra $A$ such that $N_\lambda \subseteq N_\mu$.  With the notations as in Definition \ref{Transfer-of-measures} and Remark \ref{notation-agreement}, we say
		\begin{enumerate}[(i)]
			\item   $\mu: A \longrightarrow \mathbb{C}$ is  {\it weak* continuous} with respect to  $\lambda$ if $\mu$ has a weak* continuous extension to a PLF $\hat{\mu}$ on $L^{\infty}({\lambda})$, i.e., the extension $\hat{\mu}:L^{\infty}({\lambda})\longrightarrow \mathbb{C}$ is weak*-continuous. We use $\mu \ll_{w} \la$ to show that $\mu$ is weak* continuous with respect to $\la$, and we denote the set of all weak*-continuous PLFs on $A$ with respect to $\la$ by $WC[\la]$.
			
			\item   $\mu$ is called {\it weak*-singular} with respect to $\lambda$ if no extension of $\mu$ to $L^\infty(\lambda)$ can majorize a non-zero weak*-continuous PLF.  We use $\mu \perp_{w} \la$ to show that $\mu$ is weak* singular with respect to $\la$, and we   denote the set of all weak*-singular PLFs on $A$ with respect to $\la$ by $WS[\la]$.

			\item  The zero functional is both weak* continuous  and weak* singular with respect to $\lambda$. 
		\end{enumerate}
		
	\end{defn}
	\begin{remark} \label{wk-str-continuous-wrt-itself}
		\begin{enumerate}[(i)]
			\item      Every PLF $\lambda: A \longrightarrow \mathbb{C}$ is weak* continuous with respect to itself. This is because by the GNS construction we have 
			$$\lambda^{\prime}(\pi_\lambda(a))=\lambda(a)=\langle \pi_\lambda(a)(1+N_\lambda),~ 1+N_\lambda\rangle_{L^2(\lambda)}$$
			So, by our notation in Remark \ref{notation-agreement}, we have $\lambda(a)=\langle a+N_\lambda,~ 1+N_\lambda\rangle_{L^2(\lambda)}$, which is a vector state on $L^\infty(\la)$ by \cite[Theorem 4.2.6]{Murphy}. So, by  by \cite[Theorem II.2.6]{Takesaki_I}  $\lambda$ is indeed a weak* continuous functional on $A$. 
			\item Whenever $\mu$ is weak*-continuous with respect to $\lambda$, by \cite[Theorem 7.1.12]{Kadison_Fundamentals_II} or \cite[Theorem 2.2.6]{Takesaki_I} we can find an orthogonal sequence of vectors $(\xi_n)_{n=1}^{\infty}$ in $L^2(\lambda)$ with the property $\sum_{n=1}^{\infty}  \|\xi_n\|^2 =1$ such that
			$$\hat{\mu}(a)=\sum_{n=1}^{\infty} \langle a\xi_n , ~\xi_n \rangle_\lambda ~, \quad \forall a\in L^\infty(\lambda).$$
		\end{enumerate}
		This formula is true on $\pi_\lambda(A)$. Thus, by our agreement in Definition \ref{Transfer-of-measures}, we can apply the same formula for $\mu$ on $A$. This shows the weak*-continuous functionals with respect to $\lambda$ are special among all bounded functionals.  
	\end{remark}
	Recall that a sub-cone $\mathcal{S}$  of a positive cone $\mathcal{C}$ is hereditary if whenever $s \in \mathcal{S}$ and $c \in \mathcal{C}$ with $c \leq s$, then $c \in \mathcal{S}$.
	\begin{thm} \label{properties-of-wk-ast-decomposition}
		Let $\mu$, $\lambda$ and $\tau$ be  three PLFs  on a unital C*-algebra $A$ such that $N_\lambda \subseteq N_\mu$ and $N_\lambda \subseteq N_\tau$. Then,
		\begin{enumerate}[(i)]
			\item $\mu$ has a unique weak* Lebesgue decomposition 
			$$\mu=\mu_{wc}+\mu_{ws}$$
			with respect to $\lambda$, where $\mu_{wc}$ and $\mu_{ws}$ are PLFs and $\mu_{wc}$ is the maximal weak* continuous PLF such that $\mu_{wc} \leq \mu$.
			\item $WC[\la]$ and $WS[\la]$ are hereditary cones in $A^{*}_{+}$. 
			\item $(\mu+\tau)_{wc}=\mu_{wc}+ \tau_{wc}$ and $(\mu+\tau)_{ws}=\mu_{ws}+ \tau_{ws}$
			\item  $\mu \leq  \tau$ implies $\mu_{wc} \leq  \tau_{wc}$ and $\mu_{ws} \leq  \tau_{ws}$.
		\end{enumerate}
	\end{thm}
	\bp
	\begin{enumerate}[(i)]
		\item  Let $\hat{\mu}$ be an extension of $\mu$ on $L^\infty(\lambda)$ as in Proposition \ref{wel-definedness}. By a general result in von Neumann algebra theory, see \cite[p. 127 and Theorem III.2.14]{Takesaki_I}, we have
		$$L^{\infty}(\lambda)^*= L^{1}(\lambda) \oplus_1  L^{1}(\lambda)^{\perp},$$
		in which, both summands, i.e.,  the weak* continuous $L^{1}(\lambda)$ and the weak* singular PLFs $L^{1}(\lambda)^{\perp}$, are normed closed sub-spaces of $L^{\infty}(\lambda)^*$. Hence, the decomposition $\hat{\mu}=\hat{\mu}_{wc}+\hat{\mu}_{ws}$ follows. Furthermore, since $\hat{\mu}$ is positive,  both summands are positive by \cite[Theorem 10.1.15-iii]{Kadison_Fundamentals_II}. In addition, the decomposition into weak* continuous and weak* singular parts is unique by \cite[Theorem III.2.14]{Takesaki_I}. Thus, we observe that $\hat{\mu}_{wc}$ is the maximal weak* continuous part of $\hat{\mu}$ on $L^\infty(\lambda)$, and hence also on $\pi_\lambda(A)$, as the latter is weak* dense in the former. Given that $\mu \cong \hat{\mu}|_{\pi_{\lambda}(A)}$ (see Definition \ref{Transfer-of-measures}), and that we may define $\mu_{ws} := \mu - \mu_{wc}$, it follows that the decomposition $\mu = \mu_{wc} + \mu_{ws}$ is the maximal weak* decomposition on $\pi_\lambda(A)$; hence, it is unique on it.
		
		\item   Note that by \cite[p. 127 and Theorem III.2.14]{Takesaki_I} and \cite[Theorem 10.1.15-iii]{Kadison_Fundamentals_II},  $WC[\la]$ and $WS[\la]$ are positive  cones. So, we need only to prove the hereditary properties. We first show that  $WS[\la]$ is hereditary. Let $\nu$ be a PLF and $\mu$ be a weak* singular PLF dominating $\nu$. Since $\nu=\nu_{wc}+\nu_{ws} \leq \mu=\mu_{ws}$, we see that $\nu_{wc} \leq \mu$. However, $\mu$ is singular, so by the definition of a singular PLF, we must have  $\nu_{wc}=0$. Note that we can prove this in another way. In fact, since $\mu$ is singular, by \cite[Theorem III.3.8]{Takesaki_I} we see that for any non-zero projection $e \in L^{\infty}(\la)$ there is a non-zero projection $e_0 \leq e$ such that $\mu(e_0)=0$. Hence, $\nu_{wc}(e_0)=0$, which makes $\nu_{wc}$ into a singular PLF. This is only true when $\nu_{wc}=0$. Thus, $\nu=\nu_{ws}$, so $\nu \in WS[\la]$. 
		
		Now, let $\nu$ be a PLF and $\mu$ be a weak* continuous PLF dominating $\nu$. Since $\nu=\nu_{wc}+\nu_{ws} \leq \mu=\mu_{wc}$, we see that $0\leq \nu_{ws} \leq \mu_{wc}$. Hence, $\ga=\mu_{wc}-\nu_{ws}$ is a  PLF. By the uniqueness of the decomposition into weak* continuous and weak* singular parts, see \cite[p. 127 and Theorem III.2.14]{Takesaki_I}, we obtain that $\ga_{wc}=\mu_{wc}$ and $\ga_{ws}=-\nu_{ws}$. However, this means that the singular part of $\ga$ is negative, which contradicts \cite[Theorem 10.1.15-iii]{Kadison_Fundamentals_II}. 
		
		\item Note that $\mu+\tau=(\mu+\tau)_{wc}+(\mu+\tau)_{ws}$. Also, $\mu+\tau=(\mu_{wc}+ \mu_{ws})+ (\tau_{wc}+\tau_{ws})$. So,
		$$\mu+\nu=(\mu+\tau)_{wc}+(\mu+\tau)_{ws}=(\mu_{wc}+ \tau_{wc})+ (\mu_{ws}+\tau_{ws})$$
		However, note that the weak* continous and singular bounded functionals form normed closed sub-spaces of $L^{\infty}(\lambda)^*$. Hence, by uniqueness of the decomposition the result follows. Here, the hereditary properties can also be used to show the equality. 
		\item For the last assertions note that 
		$$\mu_{wc} \leq \mu_{wc}+ \mu_{ws} \leq \tau_{wc}+\tau_{ws}$$
		and
		$$\mu_{ws} \leq \mu_{wc}+ \mu_{ws} \leq \tau_{wc}+\tau_{ws}.$$
		So, $(\tau_{wc}-\mu_{wc})+\tau_{ws}$ and $\tau_{wc}+(\tau_{ws}-\mu_{ws})$  are PLFs, with unique decompositions into positive weak* continuous and positive weak* singular PLFs, see \cite[p. 127 and Theorem III.2.14]{Takesaki_I} and \cite[Theorem 10.1.15-iii]{Kadison_Fundamentals_II}. Thus, $\tau_{wc}-\mu_{wc} \geq 0$ and $\tau_{ws}-\mu_{ws}\geq 0$.
	\end{enumerate}
	\ep
	Now, we want to show that our weak* decomposition method of PLFs is equivalent  to other known decomposition theories like   Arveson-Gheondea-Kavruk Lebesgue (AGKL) \cite{GK-ncld}. To do so, we noticed that our splitting PLF $\lambda$ should have some kind of weak tracial property known as the Kubo-Martin-Schwinger (KMS) condition defined below. At the moment, we are unable to relax this condition in our theory. We point out that in the classical setting  we are working with abelian C*-algebras, so the KMS condition is automatic.  On the other hand, a splitting PLF with the KMS property provides some sharp equations like $(\mu+\nu)_{ac}=\mu_{ac}+\nu_{ac}$, which cannot be achieved within the AGKL theory, see \cite[Corollary 3.7-ii]{GK-ncld}. After proving our equivalence theorem, we can compute the Radon-Nikodym theorem for weak* continuous PLFs. 
	
	For a C*-algebra $A$, an automorphism is defined as a $\ast$-homomorphism that is also a bijection, i.e., a $\ast$-isomorphism. The group of all such automorphisms of $A$ is denoted by $\mathrm{Aut}(A)$. In the context of von Neumann algebras, automorphisms are further required to be normal, meaning they preserve the weak*-topology topology. 
	A strongly continuous one parameter group homomorphism  is a mapping $\sigma: \mathbb{R} \longrightarrow \mathrm{Aut}(A)$ which is continuous with respect to the Euclidean topology of $\mathbb{R}$ and the strong operator topology topology of $\mathrm{Aut}(A)$. If we put $\sigma_t=\sigma(t)$, the continuity of $\sigma$ means that for any $a \in A$ the mapping $ \mathbb{R} \longrightarrow A$ given by $t \longrightarrow \sigma_t(a) $ is Euclidean-norm  continuous. For von Neumann algebras, we require the continuity of the latter mapping when $A$ has the strong operator topology. In this situation, 
	$\sigma$ is called the time evolution of the   C*-dynamical (or the W*-dynamical) system  $(A, \mathbb{R} , \sigma) $, see also the Remark \ref{C-ast-dyn-sys} and Definition \ref{W-ast-cross-prd} for the definitions of the latter systems.
	
	The following definition is based on \cite[Definition 5.3.1]{Bratteli-OA-II}   for the KMS states in quantum mechanics.
	
	\begin{defn} \label{KMS-state}
		Let    $(A, \mathbb{R} , \sigma) $ be a C*-dynamical system. A PLF $\la$ on $A$ is called  a Kubo-Martin-Schwinger (KMS) PLF for $(A, \mathbb{R} , \sigma) $ if 
		\begin{enumerate}[(i)]
			\item for each $a \in A$ the mapping $ \mathbb{R} \longrightarrow A$ given by $t \longrightarrow \sigma_t(a) $ has an analytic extension to a mapping $f_a:  \{z \in \mathbb{C}~:~ 0 < \mathrm{Im}(z) < 1 \}  \longrightarrow A$, and define $\sigma_z(a):=f_a(z) $,
			\item the weak tracial property $\la(ab)=\la(b\sigma_{i}(a))$ holds.
			\item $\la$ is time invariant, i.e.,  $\lambda(\sigma_t)=\lambda$ for any $t \in \mathbb{R}$
		\end{enumerate}
		We use the acronym KMS-PLF for  a Kubo-Martin-Schwinger positive linear functional. 
	\end{defn}
	\begin{remark}
		Usually one extracts  the dynamics $\sigma$ from the PLF $\lambda$.    Note that according to \cite[Proposition 5.3.3]{Bratteli-OA-II}, condition (i) and (ii) imply condition (iii). Also, by \cite[Chapter 5]{Bratteli-OA-II} a weaker form of condition (i) is sufficient for the definition of a KMS-PLF, i.e., one requires that the members of a (weak*) dense $\ast$-sub-algebra have analytic extensions to the strip $\{z \in \mathbb{C}~:~ 0 < \mathrm{Im}(z) < \beta \}  $, see \cite[P. 77]{Bratteli-OA-II}  and  \cite[P. 75]{Naaijkens-Qspin}. In this situation, we should replace (ii) by $\la(ab)=\la(b\sigma_{\beta i}(a))$.  The KMS condition puts a restriction on the non-commutativity of the PLF $\lambda$ without requiring it to be a full trace. This is important since most C*-algebras do not have a trace, but they can posses a KMS state. Here, $\beta$ measures how far $\lambda$ is being from an actual trace. For $\beta =0$, we have an actual trace. 
	\end{remark}
	
	Before giving the next example, we need some definitions and constructions related to C*-algebras.
	The inductive or direct limit of \(C^*\)-algebras is a construction used to build larger \(C^*\)-algebras from a sequence of smaller \(C^*\)-algebras, see \cite[Chapter 6]{Murphy}. 
	
	\begin{defn} \label{inductive-limit-C-ast-alg}
		Suppose that we have a sequence of C*-algebras $(A_n)_{ n\in \mathbb{N}}$, and assume that when $n \leq m$ there is $\ast$-homomorphism \(\phi_{n,m}: A_m \to A_n\). Define
		$$
		A^{\prime} = \{ (a_n) \in \prod_{n\in \mathbb{N}} A_n ~:~ \exists ~ K \in \mathbb{N}~ \forall n\geq m \geq K,~ a_n = \phi_{n,m}(a_m) \}
		$$
		Since $\ast$-homomorphisms between C*-algebras are norm decreasing, for  any  $a=(a_n)_{n\in \mathbb{N}} \in A^{\prime}$ the sequence of norms $(\|a_n\|)_{n\in \mathbb{N}}$ is eventually decreasing.  Hence, the limit 
		$$
		p (a)  = \lim_{n \to \infty}\|a_n\|
		$$
		exists, and defines a C*- semi-norm on $A^{\prime}$. Now, we can divide $A^{\prime}$ by the kernel of $p$, and then complete it to a C*-algebra $A$. This $A$  is called the {\it inductive} or {\it direct limit} of $(A_n)_{ n\in \mathbb{N}}$ and is denoted by  \(\varinjlim A_n\). One can adjust this definition for a net of C*-algebras with appropriate *-homomorphisms among them.   
	\end{defn}

	The following example is taken from \cite[P. 13]{Bratteli-OA-I} and  \cite[P. 77]{Naaijkens-Qspin}, which is important in quantum spin systems.
	\begin{eg}
		Consider the lattice $\mathbb{Z}^n$, and to  each point of $x \in \mathbb{Z}^n$ choose a copy of a finite dimensional Hilbert space $\mathcal{H}_x$. For each finite subset $\mathcal{F}$ of $\mathbb{Z}^n$, we form the tensor product $\mathcal{H}_{\mathcal{F}}=\otimes_{x \in \mathcal{F}} \mathcal{H}_x  $ and the operator algebra  $A(\mathcal{F})=\mathscr{B}(\mathcal{H}_{\mathcal{F}})$. The direct limit of $A(\mathcal{F})$, denoted by $\mathcal{A}$, is  called the C*-algebra of local observable, \cite[P. 56]{Naaijkens-Qspin}. Note that $\mathcal{A}$ is a sub-algebra of a suitable $\mathscr{B}(\mathcal{H})$, where $\mathcal{H}$ is related to the infinite tensor product of Hilbert spaces.  One can show that there is a closed unbounded operator $T$, known as the Hamiltonian of the system, such that $e^{-\beta T}$ is a trace class operator and $\lambda(a)=\mathrm{tr}(ae^{-\beta T})$ is a KMS-PLF on $\mathcal{A}$. Note that the domain $\mathcal{D}$ of $T$  is dense in the Hilbert space $\mathcal{H}$, and $e^{-\beta T}$ is only invertible on $\mathcal{D}$.  Thus, $e^{-\beta T}$ being a trace class operator does not contradict the fact that the ideal of trace class operators contains no invertible everywhere defined operator, see \cite[P.85]{Bratteli-OA-I}. Such an operator is called a {\it density matrix} for $\la$, see \cite[P. 76]{Bratteli-OA-I}.  Also, $\sigma_t(a)=e^{itT}a e^{-itT}$ determines the dynamics of the system, see \cite[P.76]{Bratteli-OA-II}. 
	\end{eg}
	\begin{eg}
		Let $\mathscr{M}$ be a von Neumann algebra inside some $\mathscr{B}(\mathscr{H})$, and suppose that $\mathscr{M}$ has a separating and cyclic unit vector  $\eta$. Then, by \cite[Chapter VI]{Takesaki_II} or \cite[P. 285]{Bratteli-OA-I}, one can show that the vector state $\phi(a)=\langle a\eta, ~ \eta\rangle$ is a KMS state. In this situation, we let $S$ be the closure of the densely defined operator 
		\begin{eqnarray*}
			S_0 : \mathscr{M}\eta &\longrightarrow &\mathscr{H}\\
			m\eta &\longrightarrow & m^*\eta,
		\end{eqnarray*}
		and we put $\nabla=S^*S$. By the Tomita-Takesaki theory \cite{Takesaki_II},  the time evolution is given by (see \cite[P. 284]{Bratteli-OA-I})
		\begin{eqnarray*}
			\sigma: \mathbb{R} &\longrightarrow &\mathrm{Aut}(\mathscr{M}) \\
			\sigma_t(a)&=&\nabla^{it}a\nabla^{-it} ~.\blacksquare
		\end{eqnarray*}
		
	\end{eg}
	Let's recall some other forms of absolute continuity and singularity for PLFs. In the classical measure theory, absolute continuity of measures on $\mathbb{T}$ can be translated into the language of  convergence of measures. Namely, given two finite positive regular Radon measures $\mu$ and $\lambda$ on $\mathbb{T}$, one can show that $\mu$ is absolutely continuous (AC) with respect to $\lambda$, in notations $\mu \ll \la$, if there is an increasing sequence of finite positive regular Radon measures $\mu_{n}$, which are dominated by $\lambda$, and increasing monotonically to $\mu$:
	\begin{eqnarray*}
		0 \leq \mu_{n} \leq \mu, \quad \mu_{n} \uparrow \mu, \\
		\forall~n~\exists~  t_{n} >0 , \quad  \mu_{n} \leq t_{n} \lambda.
	\end{eqnarray*}
	
	The following definition, provided in \cite{GK-ncld}, generalizes the above concept to the non-commutative setting. 
	
	\begin{defn}\label{AC-msr} Consider two PLFs $\mu$ and $\lambda$ on a unital C*-algebra $A$ such that $N_\lambda \subseteq N_\mu$. 
		\begin{enumerate}[(i)]
			\item   We say $\mu$ is  {\it absolutely continuous} with respect to  $\lambda$ in the sense of Gheondea-Kavruk  (GK), denoted by $\mu \ll \la$,  if 
			\begin{itemize}
				\item there is an increasing sequence $(\mu_n)_{n\in \mathbb{N}}$ of PLFs such that $\mu_n \uparrow \mu$ point-wise, i.e., in the weak* topology of $A^{*}$ 
				\item for each $n\in \mathbb{N}$, there is  $t_n >0$ such that $\mu_n \leq t_n \lambda$
			\end{itemize}
			
			We denote the set of all absolutely continuous PLFs with respect to $\la$ by $AC[\la]$.
			\item  $\mu$ is called {\it singular} with respect to  $\la $, denoted by $\mu \perp \la$, if the only PLF majorized by both $\mu$ and $\la$ is the zero PLF. We   denote the set of all singular PLFs with respect to $\la$ by $SG[\la]$. 
		\end{enumerate}
	\end{defn}
	
	\begin{defn}
		A {\it densely defined operator} $T$ on a Hilbert space $\mathscr{H}$  is a linear operator $T: \mathscr{D}(T) \longrightarrow \mathscr{H} $ where its domain $\mathscr{D}(T)$ is a dense linear subspace of $\mathscr{H}$. We always assume such an operators is unbounded; otherwise, it possesses a bounded extension to the whole space $\mathscr{H}$, and its theory reduces to that of bounded operators.  Since densely defined unbounded operators depend so much on their specific domains, the algebraic operations among  them follow the same rules of ordinary functions. For example, if $S: \mathscr{D}(S) \longrightarrow \mathscr{H} $ is another densely defined operator, then 
		\begin{eqnarray*}
			\mathscr{D}(T+S)&=&\mathscr{D}(T) \cap \mathscr{D}(S) \\
			\mathscr{D}(TS)&=& \{~ \xi \in \mathscr{D}(S)~:~S\xi \in \mathscr{D}(T) ~\}.
		\end{eqnarray*}
		However, we note that these new sets might not be dense in $\mathscr{H}$. We say that $S$ is an extension of $T$ if $\mathscr{D}(T) \subseteq \mathscr{D}(S)$ and $S=T$ on $\mathscr{D}(T)$. The operator $S$ is {\it closed} if its graph
		\begin{equation*}
			\Gamma(S) = \{~\xi \oplus S\xi ~ : ~ \xi \in \mathscr{D}(S) ~ \}
		\end{equation*}
		is closed in $ \mathscr{H} \oplus \mathscr{H}$. The operator $T$ is {\it closable} if it has a closed extension $S$. The {\it adjoint} $T^*$ of $T$ is specified by its domain
		\begin{equation*}
			\mathscr{D}(T^*) = \{~\xi \in \mathscr{H} ~ : ~\eta \longrightarrow \langle T\eta , \xi \rangle~ \text{is a bounded linear functional on} ~\mathscr{D}(T)~ \}.
		\end{equation*}
		When $\xi \in \mathscr{D}(T^*)$, the relation $\langle T\eta , \xi \rangle = \langle \eta , T^* \xi \rangle$ holds. One can show that $T^*$ is always closed.   The operator $T $ is {\it self-adjoint} if $T=T^*$ and $\mathscr{D}(T)=\mathscr{D}(T^*)$.  The {\it resolvent set} $\rho(T)$ for $T$ is the set of all $\la \in \mathbb{C}$ such that $\la I - T$ has  a bounded inverse. The spectrum of $T$ is the set $\sigma(T)=\mathbb{C}\setminus\rho(T)$. For self-adjoint operators the spectrum contains only real numbers. There is a version of spectral theorem and functional calculus for self-adjoint unbounded operators, and based on that one can prove a polar decomposition for closed unbounded operators. A self-adjoint operator $T$ is {\it positive} if  $\langle T\xi , \xi \rangle \geq 0$ for any $\xi \in \mathscr{D}(T)$. When $T$ is positive, we have $T=\sqrt{T} \circ \sqrt{T}$, so $\mathscr{D}(T) \subseteq \mathscr{D}(\sqrt{T})$.$\blacksquare$
	\end{defn}
	
	Unbounded closed operators arise in the theory of derivatives and differential operators. One can use the result of \cite{GK-ncld} to find the Radon-Nikodym derivative for the weak* continuous part of the weak*-Lebesgue decomposition of $\mu$ with respect to $\la$.  Let $D$ be a densely defined unbounded operator on a Hilbert space $\mathscr{H}$, and $\mathscr{M}$ be a von Neumann algebra inside $\mathscr{B}(\mathscr{H})$. We say that $D$ is affiliated with $\mathscr{M}$, and we write $D \sim \mathscr{M} $ if for any $B \in \mathscr{M}^{\prime}$, we have $BD \subseteq DB$. The latter means $\mathrm{Dom}(D)$ is $B$-invariant, and $BDh=DBh$ for any $h \in \mathrm{Dom}(D)$. What is important here is that $D$ commutes with elements of $\mathscr{M}^{\prime}$ over the domain of $D$. When, $D$ is bounded, the condition $D \sim \mathscr{M} $ means $D \in \mathscr{M}$.
	
	The following  is \cite[Theorem 2.11 ]{GK-ncld}.
	\begin{thm} \label{AC-derivative}
		Let    $A $ be a unital C*-algebra. Assume that $\mu$ and $\la$ are two  PLFs on it and $\mu \ll \la$. Then, there is a unique and  possibly unbounded operator $D_{\la} \mu$ on $L^2(\la)$ such that
		\begin{enumerate}[(i)]
			\item $D_{\la} \mu$ is positive.
			\item  $D_{\la} \mu$ is affiliated with the von Neumann algebra $L^\infty(\la)^{\prime}$.
			
			\item $\pi_{\la}(A) \xi_\la$ is a subset of the domain of $\sqrt{D_{\la} \mu}$ . Here, $\xi_\la=1+N_\la$ is the cyclic vector of the GNS representation of $\la$.
			\item $\mu(a)= \langle \pi_\la(a) \sqrt{D_{\la} \mu}~\xi_\la ~, ~ \sqrt{D_{\la} \mu}~\xi_\la \rangle_{\la}  $.
		\end{enumerate}

	\end{thm}
	Here, we should note that \cite{GK-ncld} works with more general positive maps, i.e., the ones whose ranges are the full operator algebra $\mathscr{B}(\mathscr{H})$ for some Hilbert space $\mathscr{H}$, and one has to use the Stinespring's construction instead of the GNS one. However, one can simplifythe notations and formulas of \cite[Theorem 2.11]{GK-ncld}  according to \cite[Chapter 4, P.45]{Paulsen-cb}.
	
	\begin{thm} \label{equivalence-wk-str-GK}
		Let    $(A, \mathbb{R} , \sigma) $ be a C*-dynamical system. Assume that $\mu$ and $\la$ are two PLFs on it with, $N_\la \subseteq  N_\mu$. Suppose further that $\lambda$ is a KMS-PLF. Then, 
		\begin{enumerate}[(i)]
			\item $\mu \ll_{w} \la$ if and only if  $\mu \ll \la$, and
			\item $\mu \perp_{w} \la$ if and only if $\mu \perp \la$.
		\end{enumerate}
	\end{thm}
	\bp
	\begin{enumerate}[(i)]
		\item Suppose that  $\mu \ll_{w} \la$. By Definition \ref{wc-msr}, $\mu$ has a weak* continuous extension to a PLF on $L^{\infty}({\lambda})$. By \cite[P. 156]{Takesaki_I}, there is a sequence $(\xi_n)_{n\in \mathbb{N}}$ of  vectors in $L^2(\la)$ such that $s=\sum_{n=1}^{\infty} \| \xi_n\|^2 < \infty$ and $\mu(y)=\sum_{k=1}^{\infty} \langle y \xi_k, ~ \xi_k \rangle_\la $. In Remark \ref{wk-str-continuous-wrt-itself}-ii, we mentioned that we can apply this formula to weak* continuous functionals on $A$. Define a sequence of linear functionals on $A$ by
		\begin{eqnarray*}
			\mu_n(y)=\sum_{k=1}^{n}  \langle y \xi_k, ~ \xi_k \rangle_\la, \quad \forall ~ n\in \mathbb{N}
		\end{eqnarray*}
		We show that $\mu_n$'s are the desired PLFs, which satisfy the conditions of absolute continuity in Definition \ref{Arens-products}. It is obvious from the construction  that $\mu_n \uparrow \mu$. So, for any $n\in \mathbb{N}$, we only need to find a $t_n >0$ such that $\mu_n \leq t_n \la$ on positive elements like $y=a^*a$. 
		\begin{itemize}
			\item   Suppose that  $\la(y)=0$. So, $a \in N_\la \subseteq N_\mu$; thus, $a\in N_\mu$ and $\mu(y)=\mu(a^*a)=0$. So, $\mu_n(y)=0$, and the inequality $\mu_n \leq t_n \la$ holds for any positive $t_n$. 
			
			\item For the case $\lambda(y)>0$, define a bi-linear form $V$ on $A$  by 
			\begin{eqnarray*}
				V: A \times A &\longrightarrow& \mathbb{C}\\
				V(a,b)&=&\frac{\la(bya)}{\la(y)}
			\end{eqnarray*}
			Recall that $\la$ is a KMS-PLF, so for $u=b$ and $v=ya$, we have $\la(uv)=\la(v\sigma_i(u))$; hence,
			$$V(a,b)=\frac{\la(bya)}{\la(y)}=\frac{\la(ya \sigma_i(b))}{\la(y)}=\frac{\la \cdot y(a \sigma_i(b))}{\la \cdot y(1)}$$
			Since $y$ is positive, the functional $\la \cdot y$ is positive, and by Proposition \ref{norm-of-pstv-lin-fun}, 
			$$\| \la \cdot y\|=\la \cdot y(1)=\la(y).$$
			Therefore,
			$$|V(a,b)| \leq \frac{\|\la \cdot y\| \|a \sigma_i(b))\| }{|\la(y)|}= \|a \sigma_i(b))\| \leq \|a\|\|b\|$$
			Consequently, $V$ is a bounded bi-linear form on the C*-algebra $A$. By Generalized Grothendieck inequality \cite{Haag-Groth.inq}, we can find states $\phi_i$, ~$i=1, \cdots, 4$ , on $A$ such that
			$$|V(a,b)| \leq [\phi_1(a^*a)+\phi_2(aa^*)]^{\frac{1}{2}}[\phi_3(b^*b)+\phi_4(bb^*)]^{\frac{1}{2}}$$
			Applying this to $V(x_k, x_k^*)=\frac{\la(x_k^*yx_k)}{\la(y)}$, we get
			$$|\frac{\la(x_k^*yx_k)}{\la(y)}| \leq [\phi_1(x_k^*x_k)+\phi_2(x_kx_k^*)]^{\frac{1}{2}}[\phi_3(x_kx_k^*)+\phi_4(x_k^*x_k)]^{\frac{1}{2}} \leq 2\|x_k\|^2$$
			Thus,
			\begin{equation} \label{good-CBS-like-inq}
				|\la(x_k^*yx_k)| \leq 2\|x_k\|^2 \la(y).
			\end{equation}
			In addition, for any $\xi_k \in L^2(\la)$, there is a sequence $(x_{k,j})_{j=1}^{\infty}$ such that $\xi_k = \lim_{j \to \infty} x_{k,j} +N_{\la}$ in the Hilbert space $(L^2(\la) , \| \cdot\|_\la)$. So, the norm of sequence $(x_{k,j})_j$ has a positive bound $s_k \in \mathbb{R}$. However, 
			\begin{eqnarray*}
				\mu_n(y)&=&\sum_{k=1}^{n} \langle y \xi_k, ~ \xi_k \rangle_\la \\
				&=&  \sum_{k=1}^{n} \lim_{j \to \infty} \langle y x_{k,j} +N_\la, ~  x_{k,j} +N_\la \rangle_\la \\
				&=& \lim_{j \to \infty} \sum_{k=1}^{n} \la((x_{k,j})^{*}y x_{k,j})\\
				&\leq & 2\lim_{j \to \infty} \left(\sum_{k=1}^{n} \| x_{k,j}\|^2\right) \la(y),  \quad \text{by Inequality~} \ref{good-CBS-like-inq}\\
				&\leq & 2 \left(\sum_{k=1}^{n} s_{k}^{2}\right) \la(y)\\
				&\leq & t_n \la(y).
			\end{eqnarray*}
			
			So, by letting $t_n = 2 \left(\sum_{k=1}^{n} s_{k}^{2}\right) $, we see that $\mu_n \leq t_n \la$. Hence, the second condition of the definition of absolute continuity is satisfied. Thus, $\mu \ll \la$. 
		\end{itemize}
		On the other hand, assume that $\mu \ll \la$. By the Theorem \ref{AC-derivative}, there is a positive unbounded operator $D_{\la} \mu$ on $L^2(\la)$ affiliated with the von Neumann algebra $L^\infty(\la)^{\prime}$ such that
		$\mu(a)= \langle \pi_\la(a) \sqrt{D_{\la} \mu}~(1+N_\la) ~, ~ \sqrt{D_{\la} \mu}~(1+N_\la) \rangle_{\la}  $. By letting $\eta_\mu :=\sqrt{D_{\la} \mu}~(1+N_\la) \in L^2(\lambda)$, we have 
		$\mu(a)= \langle \pi_\la(a) \eta_\mu ~,\eta_\mu \rangle_{\la}  $. We can extend this to $L^\infty(\la)$  in an obvious way, i.e., $\hat{\mu}(y)= \langle y \eta_\mu ~,\eta_\mu \rangle_{\la}  $  By \cite[Theorem 7.1.12]{Kadison_Fundamentals_II}, $\hat{\mu}$ defines a weak*-continuous functional on $L^\infty(\la)$. Thus, $\mu \ll_{w} \la$. 
		
		\item Let $\mu \perp_{w} \la$, but $\mu \not\perp \la$. By \cite{GK-ncld}, this means that we can decompose $\mu$ against $\la$, i.e., $\mu =\mu_{ac}+\mu_s$, where $\mu_{ac}\neq 0$ and  $\mu_{ac} \ll \la$.  By part (i), we conclude that $\mu_{ac}$ is weak* continuous with respect to $\la$. Since  $\mu_{ac} \leq \mu$, we see that $\mu$ dominates the non-zero weak* continuous PLF $\mu_{ac}$. However, this is impossible as $\mu$ is weak* singular.  
		
		The reverse implication is proved similarly. Indeed, if  $\mu \perp \la$ but  $\mu \not\perp_{w} \la$, then we have a weak* decomposition $\mu =\mu_{wc}+\mu_{ws}$, where $\mu_{wc} \neq 0$. Now, $\mu_{wc} \leq \mu$ and $\mu_{wc} \ll \la$ by part (i). So, from  Definition \ref{AC-msr} we can find a non-zero PLF $\rho$ such that  $\rho \leq \mu$ and $\rho \leq \la$. However, this contradicts the assumption that $\mu \perp \la$.
	\end{enumerate}
	
	\ep
	\begin{remark}
		Let $\la$ be a faithful PLF with the tracial property $\la(xx^*)=\la(x^*x)$ on $L^\infty(\la)$.  By \cite[Lemma VII-2.1]{Takesaki_II} $\mathscr{N}_{\lambda} =\{  x \in \mathscr{M} ~:~ \lambda(x^* x) < \infty \}$ is a left ideal in $L^\infty(\la)$ and $\mathscr{M}_{\lambda} =\mathscr{N}_{\lambda}^{*}\mathscr{N}_{\lambda} =\{ \sum_{j=1}^{n} y_{j}^{*}x_j ~:~ n\in \mathbb{N}, ~ x_j, y_j \in \mathscr{N}_\lambda \}=\mathscr{N}_{\lambda}^{*}\cap \mathscr{N}_{\lambda}$ is a hereditary $\ast$-subalgebra of  $L^\infty(\la)$. By \cite[P. 320]{Takesaki_I}, $L^1( \la) $ is the completion of $\mathscr{M}_\la$ under $\|x\|_1=\la(|x|)$, and $L^\infty( \la) $ becomes an $L^1( \la) $-module. By \cite[P. 321]{Takesaki_I}, for any $\mu \ll_{wc} \la $, there is an $x \in L^1( \la)$ such that $\mu(a)=\la(ax)$. In this situation, the proof of Theorem \ref{equivalence-wk-str-GK} is similar to Theorem \ref{wc-ac-equiv--classical-measure}. 
	\end{remark}
	
	We are  now in a position to strengthen some of the results of \cite{GK-ncld}  for splitting KMS-PLFs. 
	
	\begin{thm} 
		Let    $(A, \mathbb{R} , \sigma) $ be a C*-dynamical system. Assume that $\mu$, $\nu$ and $\la$ are  PLFs on it with $N_\la \subseteq  N_\mu$ and $N_\la \subseteq  N_\nu$. Suppose further that $\lambda$ is a KMS-PLF. Then, 
		\begin{enumerate}[(i)]
			\item $AC[\la]$ is a positive hereditary cone, 
			
			\item $SG[\la]$ is a positive hereditary cone, 
			\item $(\mu+\nu)_{ac}=\mu_{ac} + \nu_{ac}$ with respect to $\la$,
			\item $(\mu+\nu)_{s}=\mu_{s} + \nu_{s}$ with respect to $\la$,
			\item If $\mu \leq  \nu$, then $\mu_{ac} \leq \nu_{ac}$ and $\mu_{s} \leq \nu_{s}$
			\item The weak* Lebesgue decomposition $\mu=\mu_{wc}+\mu_{ws}$ with respect to $\la$ and the Arveson-Gheondea-Kavruk Lebesgue (AGKL) decomposition $\mu=\mu_{ac}+\mu_{s}$ with respect to $\la$ are the same.
		\end{enumerate}
	\end{thm}
	\bp
	Combine Theorem \ref{properties-of-wk-ast-decomposition} and Theorem \ref{equivalence-wk-str-GK}.
	\ep
	The materials in the following definition can be found in \cite[Chapter X]{Conway}.

	We can now give a version of Theorem \ref{AC-derivative} for weak* continuous functionals:
	
	\begin{thm} \label{Radon-Nikodym-derivative}
		Let    $(A, \mathbb{R} , \sigma) $ be a C*-dynamical system. Assume that $\mu$ and $\la$ are  PLFs on it with $N_\la \subseteq  N_\mu$,  and $\lambda$ is a KMS-PLF. If $\mu \ll_{wc} \la$, then there is a unique and  possibly unbounded operator $D_{\la} \mu$ on $L^2(\la)$ such that
		\begin{enumerate}[(i)]
			\item $D_{\la} \mu$ is positive.
			\item  $D_{\la} \mu$ is affiliated with the von Neumann algebra $L^\infty(\la)^{\prime}$.
			
			\item $\pi_{\la}(A) \xi_\la$ is a subset of the domain of $\sqrt{D_{\la} \mu}$ . Here, $\xi_\la=1+N_\la$ is the cyclic vector of the GNS representation of $\la$.
			\item $\mu(a)= \langle \pi_\la(a) \sqrt{D_{\la} \mu}~\xi_\la ~, ~ \sqrt{D_{\la} \mu}~\xi_\la \rangle_{\la}  $.
			\item If $\mu \leq t \la$ for a positive number $t \in \mathbb{R}$, then $D_{\la} \mu$ is a bounded positive operator which belongs to $L^\infty(\la)^{\prime}$ and $\mu(a)= \langle \pi_\la(a) D_{\la} \mu (\xi_\la)~, ~ \xi_\la \rangle_{\la}  $. Also, the assignment $\mu \longrightarrow D_{\la} \mu $ is an affine mapping with respect to $\mu$.
		\end{enumerate}

	\end{thm}
	\bp
	The proofs of $(i)$ to $(iv)$ is a result of the equivalent Theorem \ref{equivalence-wk-str-GK} and \cite[ 2.11]{GK-ncld}. The final part is a result of the equivalent Theorem \ref{equivalence-wk-str-GK} and \cite[ Corollary 2.3]{GK-ncld}.
	\ep
	
	\begin{remark}
		When $\mu \ll_{wc} \la$, by Theorem \ref{equivalence-wk-str-GK} and Definition \ref{AC-msr}, we know that there is an increasing sequence $(\mu_n)_{n\in \mathbb{N}}$ of PLFs such that $\mu_n \uparrow \mu$ point-wise, i.e., in the weak* topology of $A^{*}$. Also,  for any $n\in \mathbb{N}$ there is a positive $t_n$ such that $\mu_n \leq t_n \lambda$. The latter condition implies that  $\mu_n \ll_{ac} \la$, and hence each  $D_{\la} \mu_n$ is  a positive bounded operator in $L^\infty(\la)^{\prime}$. By \cite[Part F]{GK-ncld}, one finds that $(I+D_{\la} \mu_n)^{-1} \xrightarrow{SOT} (I+D_{\la} \mu)^{-1}$. This convergence is described by saying that $D_{\la} \mu_n \xrightarrow{SRS} D_{\la} \mu$, where ``SRS" stand for the {\it  strong resolvent sense}. Note that the number $-1$ does not belong to the spectra  of $D_{\la} \mu_n$'s, so by \cite[Theorem 5.1.9]{AnalysisNow} the  operators $I+D_{\la} \mu_n$ and $I+D_{\la} \mu$ are invertible. In summary, though the Radon-Nikodym derivative  $D_{\la} \mu$ is generally an unbounded operator, one can still approximate it in a suitable sense by bounded Radon-Nikodym derivatives. 
	\end{remark}

	\begin{remark}
		Let $\mu$ and $\lambda$ be positive linear functionals defined on the Cuntz-Toeplitz operator system, as in \cite{JM-ncld} and \cite{NCLD-Fouad}. One can apply the method presented here to decompose $\mu$ with respect to $\lambda$. First, as in \cite{NCLD-Fouad}, $\lambda$ can be classified into two cases: Cuntz and non-Cuntz. Second, it is essential to ensure that the KMS condition for $\lambda$ holds on the ambient Cuntz C*-algebra.  In the non-Cuntz case, the decomposition method of \cite{NCLD-Fouad} generalizes the approach from \cite{JM-ncld}. In the Cuntz case, $\lambda$ has a unique extension to the ambient Cuntz C*-algebra by \cite[Proposition 5]{JMT-NCFM}. This allows the decomposition to be carried out in that setting, yielding the desired result. While there are technical subtleties that must be addressed, the key missing component in \cite{JM-ncld} and \cite{NCLD-Fouad} is the KMS condition. In fact, as shown in \cite[Example 1.5.8]{DNV-free.var}, the Lebesgue expectation $m$ naturally satisfies the KMS condition. This hidden property underlies the equivalence between the reproducing kernel Hilbert space and the weak* decomposition established in \cite{JM-ncld}.
		
	\end{remark}
	
	\bibliographystyle{unsrtnat}

	\Addresses
	
\end{document}